\title{Colouring complete bipartite graphs from random lists}
\author {Michael Krivelevich\thanks{
Department of Mathematics, Raymond and Beverly Sackler Faculty of
Exact Sciences, Tel Aviv University, Tel Aviv 69978, Israel.
E-mail: krivelev@post.tau.ac.il. Research supported in part by a
USA-Israel BSF Grant and by a grant from the Israel Science
Foundation.} \and Asaf Nachmias\thanks{ Department of Mathematics,
Raymond and Beverly Sackler Faculty of Exact Sciences, Tel Aviv
University, Tel Aviv 69978, Israel. E-mail:
asafnach@post.tau.ac.il.}}

\documentclass[11pt]{article}
\usepackage{latexsym}
\usepackage{amsfonts}
\usepackage{amsmath}
\newtheorem{theo}{Theorem}[section]
\newtheorem{prop}[theo]{Proposition}
\newtheorem{lemma}[theo]{Lemma}
\newtheorem{coro}[theo]{Corollary}

\date{}
\begin{document}

\maketitle

\begin{abstract}
Let $K_{n,n}$ be the complete bipartite graph with $n$ vertices in
each side. For each vertex draw uniformly at random a list of size
$k$ from a base set $\mathcal{S}$ of size $s=s(n)$. In this paper
we estimate the asymptotic probability of the existence of a
proper colouring from the random lists for all fixed values of $k$
and growing $n$. We show that this property exhibits a sharp
threshold for $k\geq 2$ and the location of the threshold is
precisely $s(n)=2n$ for $k=2$, and approximately
$s(n)=\frac{n}{2^{k-1}\ln 2}$ for $k\geq 3$.

\end{abstract}

\section{Introduction}
Let $G$ be a simple and undirected graph. Assign to each vertex
$x$ of $G$ a set $L(x)$ of colours (positive integers). Such an
assignment $L$ of sets to vertices in $G$ is referred to as a
\em{colour scheme} \rm for $G$. An \em L-colouring \rm of $G$ is a
mapping $f$ of $V(G)$ into the set of colours such that $f(x) \in
L(x)$ for all $x\in V(G)$ and $f(x)\neq f(y)$ for each $(x,y) \in
E(G)$. If $G$ admits an $L$-colouring, then $G$ is said to be \em
L-colourable \rm. In case of $L(x)=\{1,\ldots ,k\}$ for all $x \in
V(G)$, we also use the terms \em k-colouring \rm and \em
k-colourable \rm respectively. A graph $G$ is called \em
k-choosable \rm if $G$ is $L$-colourable for every colour scheme
$L$ of $G$ satisfying $|L(x)|=k$ for all $x \in V(G)$. The \em
chromatic number $\chi(G)$ \rm (\em choice number $ch(G)$ \rm) of
$G$ is the least integer $k$ such that $G$ is $k$-colourable
($k$-choosable). The choosability concept was introduced
independently by Vizing \cite{Viz} and by Erd\H{o}s, Rubin and
Taylor \cite{ERT}.

Assign colour lists to the vertices of $G$ by choosing for each
vertex $v$ its colour list, $L(v)$, uniformly at random from all
$k$-subsets of a ground set $\mathcal{S}=\{1\ldots s\}$.
Intuitively, the larger is $s$, the more spread are the colours,
and the easier it is to colour $G$ from the chosen random lists.
The question is how large should be the value of $s=s(k,G)$ to
guarantee the almost sure choosability from random lists.

In \cite{KN} the authors solve the problem for the $d$-th power of
a cycle on $n$ vertices, denoted by $C_n^d$ (i.e. the vertices of
$C_n^d$ are those of the $n$-cycle, and two vertices are connected
by an edge if their distance along the cycle is at most $d$). In
this case the authors prove that for $d\geq k$, the threshold for
choosability occurs at $n^{1/k^2}$ (note that the threshold does
not depend on $d$), i.e. that if the $s(n)=\omega(n^{1/k^2})$
almost surely a proper choice exists and if $s(n)=o(n^{1/k^2})$
almost surely a proper choice does not exist. The threshold in
this case is coarse. That is to say that if $s(n) \sim tn^{1/k^2}$
for fixed $t$ then the probability that there exists a proper
choice tends to $\phi(t)$, where $\phi(t)$ is an increasing
positive function which tends to $0$ when $t$ tends to $0$, and
tends to $1$ when $t$ tends to infinity. This is achieved in
\cite{KN} by showing that the choosability property is controlled
by a very local property: the appearance of a clique of size $k+1$
with identical lists drawn for each vertex.

In this paper we study the problem for the complete bipartite
graph with $n$ vertices in each side, denoted by $K_{n,n}$. Let
$\mathcal{S}$ be a set of colours of size $s(n)$. For each vertex
of $K_{n,n}$ draw uniformly at random $k$ colours from
$\mathcal{S}$ to form a colour scheme $L(k,s)$. Denote by
$p(n)=p(n,k,s)$ the probability that $K_{n,n}$ is
$L(k,s)$-colourable from the random lists. For $k=1$ it is an easy
exercise to check that the choosability property has a coarse
threshold at $s(n)=n^2$. We prove the following:
\begin{theo}\label{mainthm}
For $k=2$ and any $\epsilon> 0$:

$$p(n) = \left \{
\begin{array}{ll}
o(1), & s(n) \leq (2-\epsilon)n  , \\
1-o(1), & s(n) \geq (2+\epsilon)n . \\
\end{array} \right .
$$

And for $k \geq 3$ and any $\epsilon > 0$:

$$p(n) = \left \{
\begin{array}{ll}
o(1), & s(n) \leq \frac{(1-\epsilon)n}{2^{k-1}\ln 2}, \\
1-o(1), & s(n) \geq \frac{(1+\epsilon)n}{2^{k-1}\ln 2-(k+1)\ln 2-1} . \\
\end{array} \right .
$$

\end{theo}

We shall also use the standard asymptotic notation and
assumptions. In particular we assume that the parameter $n$ is
large enough whenever necessary. For two functions $f(n)$ and
$g(n)$, we write $f=o(g)$ if $\lim_{n\rightarrow\infty} f/g=0$,
and $f=\omega(g)$ if $g=o(f)$. Also, $f=O(g)$ if there exists an
absolute constant $c>0$ such that $f(n)<cg(n)$ for all large
enough $n$; $f=\Theta(g)$ if both $f=O(g)$ and $g=O(f)$ hold; and
$f\sim g$ if $\lim_{n\rightarrow\infty} f/g=1$.


%
%

\section {The general setting}

We shall use the standard probability spaces of random
hypergraphs. Denote by $H_k(n,p)$ the space of the $k$-uniform
hypergraphs on $n$ vertices where each of the $\binom{n}{k}$ edges
is selected independently with probability $p$, and by $H_k(n,m)$
the space of all the $k$-uniform hypergraphs on $n$ vertices and
$m$ edges uniformly distributed.

It is convenient to consider the following equivalent formulation:
Given two independent $k$-uniform random hypergraphs, $H_1$ and
$H_2$, distributed as $H_k(s(n),n)$ on the same vertex set,
$[s(n)]$, then $p(n,k,s)$ is the probability that $H_1$ has a
vertex cover (a set of vertices such that every edge contains at
least one vertex from the set) disjoint from a vertex cover of
$H_2$. To see the equivalence, simply let the vertex set be
$\mathcal{S}$, the edges of $H_1$ be the lists of the left side of
$K_{n,n}$ and similarly, the edges of $H_2$ be the lists of the
right side of $K_{n,n}$. $K_{n,n}$ has a proper colouring from the
lists iff $H_1 \cup H_2$ has such disjoint covers. We colour the
graph by picking for each vertex on the left side a colour present
in the corresponding edge in $H_1$ which is also a member of
$H_1$'s cover and similarly for vertices on the right side. This
justifies the following definition:
\newline

\noindent {\bf Definition}: A pair of hypergraphs on the same
vertex set is said to admit {\em disjoint covers} if there exist
two disjoint vertex covers, one for each hypergraph.
\newline

We shall often refer to edges of $H_1$ as {\em red}, to edges of
$H_2$ as {\em blue} and to disjoint covers as a function $\sigma
:[n] \to \{red,blue\}$ such that for every edge, $e\in H_1 \cup
H_2$, there exists $v \in e$ such that $\sigma(v)$ is of the same
colour as the $e$ (note that there are many possible functions).

It is also worth noting that this problem is equivalent to the
random $k$-SAT problem (see, e.g., \cite{AP} and its references)
with the restriction that every clause must have its literals with
the same sign. To see the equivalence, simply let the variables be
the set of vertices (originally $\mathcal{S}$), let every red edge
be a clause with positive literals and let every blue edge be a
clause with negative literals. Clearly disjoint covers exists iff
the formula is satisfiable.

\noindent The following theorems clearly imply Theorem
\ref{mainthm}:

\begin{theo}\label{th1}



For every fixed $k\geq 2$, let $H_1$, $H_2$ be two random
$k$-uniform random hypergraphs distributed independently as
$H_k(n,m(n))$, denote their union by $H_{1,2}(n,m(n))$. Then the
property of having disjoint covers has a sharp threshold. That is
to say that for every $0 < C < D < 1$, and for every $\epsilon >
0$, if $e:\mathbb{N} \to \mathbb{N}$ is a function such that
$$ D > Pr [H_{1,2}(n,e(n)) \textrm{ admits disjoint covers}] > C, $$
then if $m(n) \leq (1-\epsilon) e(n)$, $H_{1,2}(n,m(n))$ admits
disjoint covers almost surely, and if $m(n) \geq (1+\epsilon)
e(n)$, $H_{1,2}(n,m(n))$ admits no disjoint covers almost surely.

\end{theo}

\begin{theo}\label{th2}
Let $G_1$, $G_2$ be two random graphs on the vertex set $[n]$,
distributed independently as $G(n,m(n))$. Then for any
$\epsilon>0$:

$$ Pr[G_1 \textrm{ and } G_2 \textrm{ admit disjoint covers}] =
\left \{
\begin{array}{ll}
1-o(1), & m(n) \leq (\frac{1}{2}-\epsilon)n, \\
o(1), & m(n) \geq (\frac{1}{2}+\epsilon)n. \\
\end{array} \right .
$$
\end{theo}

\begin{theo}\label{th3}
For every fixed $k\geq 3$, let $H_1$, $H_2$ be two random
$k$-uniform random hypergraphs distributed independently as
$H_k(n,m(n)=rn)$. Then if $r > 2^{k-1}\ln 2$, $H_1$ and $H_2$
admit no disjoint covers almost surely.

\noindent Furthermore, there exists some $C=C(k)$ such that $C>0$
and if $r < 2^{k-1}\ln 2-(k+1)\ln 2-1$ then:

$$ Pr[H_1 \textrm{ and } H_2 \textrm{ admit disjoint covers}] > C.$$
\end{theo}

\noindent  An immediate corollary of Theorem \ref{th1} and Theorem
\ref{th3} is:

\begin{coro}
For every $r < 2^{k-1}\ln 2-(k+1)\ln 2-1$, if $m(n) \leq rn$, then
$H_1$ and $H_2$ admit disjoint covers almost surely.
\end{coro}

The rest of the paper is organized as follows: In Section 3 we
prove the sharpness of threshold result (Theorem \ref{th1}). In
section 4 we discuss the problem in the case $k=2$ and prove
Theorem \ref{th2}, and in section 5 we prove Theorem \ref{th3}.

%
%

\section {Sharpness of threshold}

We follow the approach outlined in Friedgut's survey on the sharp
threshold phenomenon (see \cite{Fried2}). \newline We prove the
sharpness result assuming $H_1$ and $H_2$ are distributed
independently as $H_k(n,p)$. Our probability space is $\Omega =
\{0,1\}^N$, $N=2\binom{n}{k}$, in which the first half of the
coordinates represents the edges of $H_1$ and the second half
represents the edges of $H_2$.

A subset $A$ of $\{0,1\}^N$ is called monotone increasing if
whenever $x \in A, x' \in \{0,1\}^N, x_i \leq x'_i$ for
$i=1,\ldots ,N$ then $x' \in A$. For $0\leq p \leq 1$, define
$\mu_p$, a product measure on $\{0,1\}^N$ with weights $1-p$ at
$0$ and $p$ at 1. That is to say:

$$ \mu(x) = p^{|x|}(1-p)^{N-|x|} \textrm { where } |x|=|{1\leq i
\leq n : x_i=1}|.$$

Denote by $H(n,p)$ a random element of $\Omega$ equipped with
$\mu_p$, and let $A\subset \Omega$ be monotone increasing. For
every $n$ define $f_n(p)= Pr[H(n,p) \in A]$. For every $n$, $f_n$
can easily be seen to be a monotone increasing polynomial in $p$.
For every $\alpha \in (0,1)$ define the function $p _\alpha (n)$
as the unique value for which $f_n(p_\alpha(n))=\alpha$.

In this context we say that a property, $A \subset \{0,1\}^N$, has
a sharp threshold if for every $\alpha \in (0,1)$, every $\delta
> 0$  and every $p = p(n)$,

$$ \lim_{n \to \infty} Pr[H(n,p(n)) \in A] =
\left \{
\begin{array}{ll}
1, & p(n) \geq (1+\delta)p_{\alpha}(n), \\
0, & p(n) \leq (1-\delta)p_{\alpha}(n). \\
\end{array} \right .
$$

\noindent Since the property of not having disjoint covers is
monotone, standard calculations (see, for example, \cite{JLR})
show that this implies a sharp threshold for when the distribution
is $H_k(n,m(n))$ thus implying Theorem \ref{th1}.

The sharpness of threshold follows from the condition that for all
$\alpha \in (0,1)$ and for all $C>0$ exists $n_0$ such that for
all $n>n_0$, $p_\alpha(n)f_n'(p_\alpha(n)) > C$. Indeed, for every
$0 < \epsilon < 1-\alpha$, and for every $C>0$, there exists
$p^*(n)$ such that $p_{\alpha}(n) < p^*(n) < p_{1-\epsilon}(n)$
and,

$$1-\epsilon-\alpha = f_n(p_{1-\epsilon}(n))-f_n(p_{\alpha} (n)) =
f_n'(p^*(n))(p_{1-\epsilon}(n)-p_{\alpha} (n)) > C - \frac{C
p_{\alpha}(n)}{p_{1-\epsilon}(n)},$$

\noindent thus, $\frac{ p_{\alpha} (n)}{p_{1-\epsilon}(n)} = 1 -
o(1)$, hence for every $\delta > 0$ exists $n_0$ such for all
$n>n_0$, $(1+\delta) p_{\alpha} (n) > p_{1-\epsilon}$, and thus
$f_n((1+\delta) p_{\alpha} (n)) = 1- o(1)$, as required. The
second part of the requirement of a sharp threshold is proven
similarly.
\newline

We will assume by contradiction that there exists $\alpha \in
(0,1)$ such that $p_\alpha(n)f'(p_\alpha(n)) < C$ for some $C$.
Our proof relies on the following result of Bourgain \cite{Fried},
which provides a sharp-threshold criteria for general monotone
properties.

\begin{theo}[Bourgain]\label{bour}
Let $A \subset \{0,1\}^n$ be a monotone property, and $C>0$
constant. Assume that there exists $\alpha \in (0,1)$ such that
$\mu_p(A)=\alpha$, $p\frac{d\mu_p(A)}{dp} < C$ and $p=o(1)$, then
there is $\delta = \delta(C) > 0$ such that either:

$$\mu_p(x\in \{0,1\}^n : x \textrm{ contains } x'\in A \textrm{ of
size } |x'|\leq 10C) > \delta,$$

\noindent or there exists $x' \not \in A$ of size $|x'|\leq 10C$
such that the conditional probability satisfies:

$$\mu_p(x\in A|x' \subset x)>\alpha + \delta.$$
\end{theo}

The idea of the proof of Theorem \ref{th1} is as follows. Assuming
the threshold for disjoint covers is coarse (i.e., not sharp), we
have by Theorem \ref{bour} that there exists a fixed set of edges,
$M$ (some red and some blue), whose addition changes the
probability of disjoint covers by some positive constant,
$\delta$. Since the property of having disjoint covers is
symmetric with respect to hypergraph automorphisms, it can easily
be shown that adding a random copy of $M$ has the same effect. On
the other hand, the fact that the threshold is coarse implies that
the addition of a large number of random edges has almost no
effect on the existence of disjoint covers. We will use the Erd\H
os-Simonovits Theorem \cite{ES} to show that these two conclusions
contradict each other.
\newline

\noindent{\bf Proof of Theorem \ref{th1}.} Our probability space
is $\Omega = \{0,1\}^{2\binom{n}{k}}$, equipped with the
probability measure $\mu_p$ defined earlier. The first half of the
coordinates represents $H_1$ and the second half represents $H_2$.
The property $A \subset \Omega$ is all the pairs ($H_1$, $H_2$)
with no disjoint covers.


Assume that the theorem does not hold, then by the above
discussion, there exist constants $\alpha$ and $C$ as in Theorem
\ref{bour}. Easy first moment calculations done in Section 5 show
that if $p=p(n)$ is such that $\mu_p(A)=\alpha$ then
$p(n)=O(n^{-k+1})$. Within this range of probabilities, it is an
easy exercise to see that for every set of vertices, $B\subset
[n]$ of size $|B|\leq 10kC$, the edges spanned on these vertices
have disjoint covers almost surely (for instance, for $k=2$ this
can only be a union of trees and unicyclic components and for
$k>2$ check that almost surely for any subhypergraph exist a
vertex of degree at most $1$). Thus, the first option of Theorem
\ref{bour} does not hold, therefore the second option must hold,
i.e., there exists some $M \in \Omega - A$ of size $|M|\leq 10C$
such that $Pr[H(n,p) \in A | M\subset H(n,p) ] > \alpha+\delta$.
Since $A$ is invariant under hypergraph vertex-automorphisms, the
latter holds for any isomorphic copy $M'$ of $M$. Because of this,
if we first draw a random copy of $M$, from all possible ones on
the vertex set $[n]$, and then draw $H \in \Omega$ their union
will be in $A$ with probability $> \alpha +\delta$ (note that the
probability space here is the product space of $\Omega$ equipped
with $\mu_p$ and the space of all isomorphic copies of $M$ on the
vertex set $[n]$, uniformly distributed). Now, since
$p\frac{d\mu_p(A)}{dp} < C$, $\lim _{\epsilon \to 0}
\frac{\mu_{p+\epsilon p}(A) - \mu_p (A)}{\epsilon} < C$, implying
that there exists $\epsilon>0$ such that $\mu_{p+\epsilon p}(A) <
\alpha +\frac{\delta}{2}$. Furthermore, the usual double exposure
routine shows that $\Omega$ equipped with the probability measure
$\mu_{p+\epsilon p}$ is probabilistically isomorphic to the union
of two independent copies of hypergraphs drawn from $\Omega$: one
with probability measure $\mu_p$ and another with $\mu_{\epsilon
'p}$ for $\epsilon '=\epsilon ' (\epsilon) > 0$ satisfying
$(1-(1+\epsilon)p)=(1-p)(1-\epsilon 'p)$. Denote by $M^*$ a
uniformly random copy of $M$ drawn from all possible copies on the
vertex set $[n]$. To sum things up, we have:


\begin{equation}\label{eq2}
Pr[H(n,p) \cup M^* \in A] > \alpha + \delta
\end{equation}

\begin{equation}\label{eq3}
Pr[H(n,p) \cup H(n,\epsilon 'p) \in A] < \alpha + \frac{\delta}{2}
\end{equation}

Now, from (\ref{eq2}) and (\ref{eq3}) we have that:

$$ \sum _{H_0 \in \Omega} Pr[H(n,p)=H_0](Pr[H_0 \cup M^* \in A] -
Pr[H_0 \cup H(n,\epsilon 'p) \in A]) > \frac{\delta}{2} $$

\noindent Therefore, exists a fixed $H_0 \in \Omega$ such that:

\begin{equation}\label{eq4}
Pr[H_0 \cup M^* \in A] - Pr[H_0 \cup H(n,\epsilon 'p) \in A] >
\frac{\delta}{2}.
\end{equation}

\noindent From this obviously, $H_0$ admits disjoint covers. We
show that this leads to a contradiction. Order the vertices of $M$
arbitrarily. Denote by $t \leq 10kC$ the number of vertices
spanned by the edges of $M$. Call an ordered t-tuple of vertices
$(v_1,\ldots ,v_t)$ in $[n]$ {\em bad} if the addition of an
ordered copy of $M$ on these vertices to $H_0$ results in a
hypergraph with no disjoint covers. Clearly, by (\ref{eq4}) at
least $\frac{\delta}{2}$ fraction of the $n^t$ ordered t-tuples of
$[n]$ are bad. We now use a small variation on a theorem of Erd\H
os and Simonovits \cite{ES}, which differs only in that it deals
with ordered $t$-tuples instead of sets of size $t$. Let $T$ be a
family of ordered $t$-tuples of distinct elements of $[n]$. An
ordered $t$-tuple $(A_1,\ldots ,A_t)$ of disjoint subsets of $[n]$
is called $T${\em -complete} if for every choice of $v_i \in A_i$,
$1 \leq i \leq t$, the resulting ordered $t$-tuple, $(v_1,\ldots
,v_t)$, belongs to $T$.

\begin{theo}[Erd\H os, Simonovits]\label{OES}
For every positive integers $k$ and $t$ and $0 < \gamma \leq 1$
there exists $\gamma ' > 0$ such that for sufficiently large $n$,
if $T \subset [n]^t$ is such that $|T|>\gamma n^t$ then with
probability at least $\gamma '$ a random choice of $t$ $k$-tuples
from $[n]$ is $T$-complete.
\end{theo}

Now, we shall get a contradiction from (\ref{eq4}). Draw uniformly
at random $t$ red edges (from $H_1$), $e_1,\ldots ,e_t$, and then
another $t$ blue edges (from $H_2$), $e_{t+1},\ldots ,e_{2t}$.
Define $T$ to be the set of all bad $t$-tuples of $[n]$. By
Theorem \ref{OES}, applied with integers $2k$ and $t$, with
probability $> \gamma '$ $(e_1 \cup e_{t+1}, e_2 \cup e_{t+2},
\ldots , e_t \cup e_{2t})$ is $T$-complete. Therefore, if we draw
$w(n) \to \infty$ red and blue edges, we will get with probability
tending to $1$ at least one occurrence of such $2t$ edges such
that  $(e_1 \cup e_{t+1}, e_2 \cup e_{t+2}, \ldots , e_t \cup
e_{2t})$ is $T$-complete. Also note that almost surely
$H_{\epsilon 'p}$ will indeed have $w(n) \to \infty$ red and blue
edges drawn, and conditioning on this event, the edges are
uniformly distributed.

Therefore, all that is left to show is that such a configuration
implies the non-existence of disjoint covers. Assume otherwise,
and let $\sigma :[n]\to \{red,blue\}$ be disjoint covers of $H_0
\cup e_1 \cup \ldots \cup e_{2t}$. Denote the vertices of $M$ by
$V(M)=\{m_1,\ldots ,m_t\}$. $M$ admits disjoint covers, $\chi
:V(M) \to \{red,blue\}$. Note that for each $1 \leq i \leq t$
there exists a vertex $v_i \in e_i \cup e_{t+i}$ such that
$\sigma(v_i) = \chi(m_i)$, thus adding an ordered copy of $M$ on
$(v_1,\ldots ,v_t)$ to $H_0$ will leave the disjoint covers
$\sigma$ intact, contradicting the fact $(v_1,\ldots ,v_t)$ is a
bad $t$-tuple. Thus, it follows that if $Pr[H_0 \cup M^* \in A] >
\frac{\delta}{2}$ then $Pr[H_0 \cup H_{\epsilon 'p} \in A] =
1-o(1)$, contradicting (\ref{eq4}).

\hfill $\Box$

%
%

\section {Case $k=2$}

In this section, we will denote by $G=G_1 \cup G_2$ our coloured
random graph where the edges of $G_1$ are red and the edges of
$G_2$ are blue. Our probability space is as usual, $\Omega =
\{0,1\}^{2\binom{n}{2}}$, equipped with the probability measure
$\mu_p$ defined earlier. The first half of the coordinates
represents the red edges and the second half represents the blue
edges. We also require the following definition:
\newline

\noindent {\bf Definition}: An {\em even alternating cycle} in $G$
is a sequence of an even number of vertices, $v_1,\ldots,v_{2t}$
such that $\{v_1, v_2\}$ is a red edge, $\{v_2, v_3\}$ is a blue
edge and so on, $\{v_{2t}, v_1\}$ is a blue edge. \newline An {\em
odd alternating cycle} in $G$ is a sequence of an odd number of
vertices $v_1,\ldots,v_{2t+1}$ such that $\{v_1, v_2\}$ is a red
edge, $\{v_2, v_3\}$ is a blue edge and so on, $\{v_{2t},
v_{2t+1}\}$ is a blue edge and $\{v_{2t+1}, v_1\}$ is a red edge.
An {\em alternating path} in $G$ is a sequence of vertices
$v_1,\ldots,v_{p}$ such that $\{v_1, v_2\}$ is a red (blue) edge,
$\{v_2, v_3\}$ is a blue (red) edge and so on.
\newline

In the proof of Theorem \ref{th2} we will see that the appearance
of large alternating cycles causes the disappearance of disjoint
covers. In fact, we will find an exact structure which determines
the existence of disjoint covers in this case.
\newline

\noindent{\bf Proof of Theorem \ref{th2}.} We will prove the
theorem assuming $G_1$ and $G_2$ are distributed $G(n,p)$ where
$p(n) = \frac{m(n)}{\binom{n}{2}}$. Since the property of having
disjoint covers is monotone, standard calculations (see
\cite{JLR}) show that this implies the theorem. The theorem
consists of two parts. We first prove that if $p(n)=\frac{c}{n}$,
such that $c$ is any fixed constant satisfying $c<1$, then
disjoint covers exist almost surely. We will present two simple
proofs for this fact, an indirect proof using Theorem \ref{th1},
and a direct proof (which also leads to a polynomial algorithm for
finding such disjoint covers). We begin with the indirect proof.
\newline

An important observation is that if $G$ has no alternating cycles,
then every subgraph of $G$ has a vertex, $v$, such that all edges
containing $v$ have the same colour. Thus we can build disjoint
covers using a greedy approach: at each step choose such $v$ and
colour it with the colour of the edges containing it. Clearly this
procedure results in disjoint covers. We will show that with
probability larger than some constant, $C>0$, $G$ has no
alternating cycles and thus, by Theorem \ref{th1} we will get
Theorem \ref{th2}. One approach is to show that the number of
vertices that participate in alternating cycles is distributed
asymptotically Poisson with constant expectation, $\mu>0$. Thus,
the probability that there exist no alternating cycles is
approximately $e^{-\mu}$. We choose to use the FKG-inequality
(see, e.g., \cite{AS}) which we quote, phrasing it according to
our needs.

\begin{theo}\label{FKG}
Let $\mathcal{A} \subset \Omega$ and $\mathcal{B} \subset \Omega$
be two monotone subsets of $\Omega$. Then:

$$ Pr[\mathcal{A} \cap \mathcal{B}] \geq
Pr[\mathcal{A}]Pr[\mathcal{B}].$$
\end{theo}

First, for any $l$ denote by $X_l$ the random variable counting
the number of alternating cycles of length $l$ in $G$. For any odd
$l$, denote by $C_1,\ldots,C_m \in \Omega$ all the possible odd
alternating cycles of length $l$. For every $1 \leq i \leq m$
define $A_i = \{G\in \Omega : C_i \not \subset G\}$, i.e., $A_i
\subset \Omega$ is the set of all graphs not containing the cycle
$C_i$. Easily, since for any $l$ vertices there are $\frac{1}{2}
(l-1)!$ possible cycle orderings, and for each $2l$ possible edge
colourings to make it an odd alternating cycle,


$$ m = \binom{n}{l} l! .$$

Furthermore, observe that $\{X_l = 0\} = \bigcap _{1\leq i \leq m}
A_i$ and that $A_i$ is a decreasing monotone subset of $\Omega$
for each $i$. Thus, using Theorem \ref{FKG} $m$ times gives:

$$ Pr[X_l = 0] \geq \left [ 1- \left ( \frac{c}{n} \right )^l
\right]^m.$$

We use the easily verifiable facts that for any $x\in [0,1/2]$,
$1-x \geq e^{-4x}$, for any natural $n>l$ $\binom{n}{l} \leq
\frac{n^l}{l!}$ to get:


$$ Pr[X_l = 0] \geq e^{-4\binom{n}{l}l!\left ( \frac{c}{n}
\right )^l} \geq e^{-4c^l}.$$

\noindent Also, a similar calculation shows that also $Pr[X_l
 = 0] \geq e^{-4c^l}$ for any even $l$.
\newline

Now, observe that the subset $\{X_l = 0\}$ is a monotone
decreasing subset of $\Omega$ for every $l$. Thus using Theorem
\ref{FKG} $n$ times,


\begin{eqnarray*}
Pr[\textrm{G admits disjoint covers}] &\geq& Pr\left [\bigcap ^{n}
_{l=1}
\{X_l=0\} \right ] \geq \prod ^{n} _{l=1} Pr[X_l=0] \\
& \geq & e^{-\sum ^{n} _{l=1} 4c^l} \geq e^{-\frac{4}{1-c}} > 0 ,\\
\end{eqnarray*}

\noindent which concludes the existential proof.
\newline

For the direct proof we find an exact structure in $G$ which
determines the existence of disjoint covers. This will imply a
polynomial algorithm for finding disjoint covers for $k=2$ of
which we omit the details (for $k>2$, deciding whether two
hypergraphs have disjoint covers is NP-hard). In \cite{APL} the
authors give a necessary and sufficient condition for a 2-SAT
formula to be satisfiable. Since, as we mentioned in Section 2,
our problem is equivalent to a certain random 2-SAT problem, we
are looking for the translation of this condition to our problem.
\newline

\noindent {\bf Definition}: An {\em odd bicycle} in $G$ consists
of two disjoint odd alternating cycles, $v_1,\ldots,v_{p}$ and
$u_1,\ldots,u_{q}$ ($q,p$ odd) such that the edges $\{v_1, v_2\}$,
$\{v_{p+1}, v_1\}$ are in the same colour, and the edges $\{u_1,
u_2\}$, $\{u_{q}, u_1\}$ are in the same colour, and one
alternating path $v_1, w_1, \ldots ,w_r, u_1$ such that the edges
$\{v_1, v_2\}$, $\{v_1,w_1\}$ are of distinct colours, and the
edges $\{u_1, u_2\}$, $\{w_r, u_1\}$ are of distinct colours as
well (the path's edges are not necessarily disjoint from the
cycles' edges).

\begin{prop}\label{BICYCLE}
$G$ admits disjoint covers if and only if $G$ does not contain an
odd bicycle.
\end{prop}
\noindent{\bf Proof.}


\newcommand {\oldproof}
{ Assume that $G$ does not contain a bicycle. We follow the greedy
approach and construct disjoint covers in the following manner:
\begin{enumerate}
\item If there exists a vertex, $v$, which all edges containing it are of the
same colour, colour $v$ with that colour and remove it and any
edge containing it from $G$. Repeat this until no such $v$ exists.
\item If there exists an odd alternating cycle, $v_1,\ldots,v_{p}$ ($p$
odd) such that the edges $\{v_1, v_2\}$, $\{v_{p+1}, v_1\}$ are of
the same colour, colour $v_1$ with that colour. If there are no
odd alternating cycles, choose an even alternating cycle and
colour an arbitrary vertex of it with red. If there exists no
alternating cycles, return to step 1.
\item Continue greedily, i.e., for every red coloured vertex, $v$,
colour any neighbour of $v$ arising from a blue edge in blue and
similarly for every blue coloured vertex. Repeat until there are
no vertices left to colour, then remove all coloured vertices and
any edge containing them and return to step 2. Return FAIL if at
some point the algorithm tries to colour an already coloured red
vertex with blue, or an already coloured blue vertex with red.
\end{enumerate}

It is easy to verify that at each point where the algorithm
colours a vertex, it does not colour the two vertices of a red
(blue) edge with blue (red). Therefore, if the algorithm does not
return FAIL, it produces disjoint covers.

Assume the algorithm returns FAIL, we will show that in that case
$G$ must contain a bicycle. Then exists $v$ such that, without
loss of generality, $v$ was coloured blue and the algorithm tried
to colour it red. Since $v$ was coloured blue, exists a sequence
of vertices which led to the colouring of $v$ with blue. Denote
that sequence by $v_1,\ldots, v_t$, such that $\{v_t,v\}$ is
coloured with blue, $v_1,\ldots, v_t, v$ is an alternating path,
$v_1,\ldots, v_t, v$ were coloured in that order by the algorithm
and $v_1$ is a member of some alternating cycle (odd or even).
Since the algorithm tried to colour $v$ red afterwards, it follows
that exists a sequence of vertices which led to the colouring of
$v$ with red, but since this happens in the same step, this
sequence must originate in one of the vertices from the former
sequence. In other words, exists $1 \leq j < t$ and vertices
$u_1,\ldots, u_p$ such that $u_1,\ldots,u_p \not \in \{v_1,\ldots,
v_t,v\}$, the edge $\{u_p,v\}$ is coloured with red, $v_j, u_1,
,\ldots, u_p, v$ is an alternating path and $v_j, u_1, ,\ldots,
u_p$ were coloured in that order by the algorithm. Choose a
maximal $j$ having these properties.

Now, note that by the definition of the algorithm, the edges
$\{v_j, v_{j+1}\}$ and $\{v_j, u_1\}$ are in the same colour,
meaning that $(v_j,\ldots,v,u_p,\ldots,u_1)$ is an odd alternating
cycle. Since we proved that an odd alternating cycle exists, this
implies that the alternating cycle $v_1$ is a member of, must be
odd. Moreover, $v_1$ must be a a member of a different odd
alternating cycle than $(v_j,\ldots,v,u_p,\ldots,u_1)$ since
otherwise $v_1 = v_j$ and by the maximality of $j$,
$(v_1,\ldots,v,u_p,\ldots,u_1)$ is the only odd alternating cycle
containing both $v_1$ and $v$, contradicting the fact that the
algorithm tried to colour $v$ in two distinct colours. This
concludes the proof of Proposition \ref{BICYCLE}. $\Box$. \newline

}

\newcommand{\newproof}
{Assume the contrary and take any minimal counter example with
respect to containment, $G$. Pick an arbitrary vertex $v$ and find
disjoint covers, $\sigma$, of $G-v$. Call an alternating path,
$v_1, \ldots, v_l$, \em $\sigma$-alternating \em if
$\sigma(v_1),\ldots, \sigma(v_l)$ is an alternating sequence of
colours. Recall that a path, $v_1, \ldots, v_l$, is called \em
simple \em if the vertices are distinct.

Denote by $P_1, \ldots, P_y$ all the $\sigma$-alternating paths
beginning in blue coloured vertices which are neighbours of $v$
arising from a red edge. Assume by contradiction that there exists
no odd alternating cycle in $G$. For any $i,j$, if $u \in P_i \cap
P_j$, then for any vertex $w$ preceding $u$ in one of the paths,
$w \in P_i \cap P_j$ since otherwise an odd alternating cycle is
formed. Denote by $V_1 \subset V$ the set of vertices of the
$P_1,\ldots, P_y$, and create a new colouring $\sigma '$ which
flips the colours of the vertices in $V_1$. Assuming there are no
odd alternating cycles in $G$, this colouring will be disjoint
covers of $G-v$. To see this note that any edge containing
precisely one vertex of $V_1$ will have the second vertex coloured
by $\sigma '$ in the same colour of the edge (otherwise, the
second vertex would have been a member in some $P_i$) and there
are no edges containing two vertices of $V_1$ such that under
$\sigma$ both vertices are coloured with the colour of the edge
(otherwise, this will result in an odd alternating cycle). Thus
$\sigma '$ are disjoint covers of $G-v$ and moreover, it is easily
seen that it can be extended to $G$ by colouring $v$ in blue, a
contradiction.

The contradiction implies that there exists an odd alternating
cycle in $G$, and moreover, there is an alternating path beginning
in a red edge containing $v$ which leads to it. The same arguments
show that there also exists an odd alternating cycle in $G$ which
begins with a blue edge containing $v$. These two obviously form
an odd bicycle. Again we have reached a contradiction. $\Box$. }


It is easy to check that an odd bicycle admits no disjoint covers.
It is left to prove that if $G$ admits no disjoint covers then it
has an odd bicycle. Let $G$ be a minimal counterexample with
respect to containment. Pick an arbitrary vertex $v\in V(G)$, and
let $\sigma:V\to \{red,blue\}$ be disjoint covers of $E(G-v)$.

A simple path $P$ in $G$ starting at $v$ is called a {\em Red
Alternating Path}, or RAP for brevity, if: 1) the colour of the
edges of $P$ alternates starting from a red edge containing $v$;
2) the colour of the vertices along $P-v$ given by $\sigma$
alternates starting from a blue neighbour of $v$. Let $V_1$ be a
set of vertices of $G-v$ reachable from $v$ by a RAP. We define a
new colouring $\sigma'$ by colouring $v$ in blue and flipping the
colours of the vertices of $V_1$. By our assumption on $G$,
$\sigma'$ does not define dijsoint covers. Therefore there exists
an edge $e=(u_1,u_2)$ missed by its corresponding colour in
$\sigma'$. Obviously, at least one of the endpoints of $e$ is in
$V_1$, say, it is $u_1$. Let $P_1$ be a RAP from $v$ to $u_1$.
Observe that the colour of $e$ should coincide with $\sigma(u_1)$.
Then, unless $u_2\in P_1$, we can extend $P_1$ to a RAP
 ending at $u_2$ by adding $e$. This shows that there is a RAP $P$
containing both $u_1$ and $u_2$, implying in particular
$u_1,u_2\in V_1$.  Since $e$ is uncovered under $\sigma'$, the
colours of $u_1,u_2$ in $\sigma$ should be identical -- and
coincide with that of $e$. But then it follows that a union of $P$
and $e$ contains an odd alternating cycle with an alternating path
connecting it to $v$, where the edge of the path containing $v$ is
red.

The same argument, with red and blue interchanged, shows that $G$
should contain an alternating cycle connected to $v$ by an
alternating path whose last edge is blue. We thus reached a
contradiction. $\Box$
\newline


%

All that is left now is to show that when $c<1$ an odd bicycle
does not exist almost surely. We use easy first moment
calculations. Denote by $Y$ the r.v. counting the number of odd
bicycles. Then,

$$ E[Y] \leq \sum _{q>2,p>2,r\geq 0}
\binom{n}{q}\binom{n}{p}\binom{n}{r} p!q!r!2 \left( \frac{c}{n}
\right )^{q+p+r+1},$$

\noindent where $p,q$ indicate the numbers of vertices in the odd
alternating cycles, and $r$ indicates the size the alternating
path connecting them. Note that after choosing the cycles'
ordering, the path's ordering and for each cycle the unique vertex
which is contained in two edges of the same colour, there are
exactly two edge colourings available. Evaluating this sum gives,

\begin{eqnarray*}
E[Y] &\leq& \frac{2c}{n}\left (\sum ^{n} _{i=0}
\binom{n}{i}i!\left (\frac{c}{n} \right)^i \right )^3 \leq
\frac{2c}{n} \left ( \sum ^{n} _{i=0} n^i
\left (\frac {c}{n}\right )^i \right )^3 \\
&\leq& \frac{2c}{n} \left( \frac{1}{1-c} \right)^3 = o(1).
\end{eqnarray*}

Thus, when $c<1$ almost surely there exists no odd bicycles, and
therefore, disjoint covers exists almost surely. \newline

We now prove the second part of Theorem \ref{th2}. The crux of the
proof is to show that when $c>1$ almost surely there exists an
alternating cycle (not necessarily a simple cycle) of length
$\omega (\sqrt n)$. The theorem follows from the observation that
by the double exposure routine, this is the same as drawing first
$G_1$ and $G_2$ distributed $G(n,\frac{c'}{n})$ where $1<c'<c$ and
then adding red and blue random edges distributed
$G(n,\frac{\epsilon}{n})$ for some $\epsilon > 0$ small enough.
Now, after the first draw we have almost surely an alternating
cycle of length $\omega (\sqrt n)$. The second draw ensures us
$\omega(1)$ random chords on that cycle. It is also easy to
observe that $\omega(1)$ random chords on an alternating cycle
results almost surely in an odd bicycle, and thus no disjoint
covers exist.

To reach an optimal result, perhaps it is best to follow the
classical paper of Ajtai, Koml\'os and Szemer\'edi \cite{AKS} in
which, using branching process techniques, they analyze the DFS
algorithm on the random graph $G(n,p=\frac{c}{n})$ where $c>1$.
This analysis leads to the existence of a cycle of linear size. It
can be shown that a slight change in their proof, in particular
drawing red/blue edges alternatingly, results in an alternating
cycle also of linear size.  However, for the sake of brevity, we
will show by easier arguments the existence of an alternating
cycle of size $\Omega(\frac{n}{\ln ^8 n})$ using a theorem of
McDiarmid \cite{MCD}. The following theorem will conclude the
proof:

\begin{theo}\label{LargeCycles}
If $c>1$, $G$ contains almost surely an alternating cycle of size
$\Omega(\frac{n}{\ln ^8 n})$.
\end{theo}
\noindent{\bf Proof.}

We associate with our $G=G_1 \cup G_2$ a new auxiliary directed
bipartite graph $H(V,E(H))$ where $V=A \uplus B$, $|A|=|B|=n$,
$A=\{a_1,\ldots, a_n\}$ and $B=\{b_1,\ldots,b_n\}$. Furthermore,
$(a_i, b_j)\in E(H)$ iff $\{i,j\}$ is a red edge in $G$ and $(b_i,
a_j) \in E(H)$ iff $\{i,j\}$ is a blue edge in $G$. Observe that a
directed cycle in $H$ corresponds to an even alternating cycle in
$G$ (again, not necessarily a simple cycle). Also, note that every
ordered pair $(a_i,b_j)$ or $(b_i, a_j)$, for $i \neq j$, is drawn
independently with probability $p$.

Denote by $B(n,n,p)$ the space of random bipartite graphs, with
$n$ vertices in each side, such that each edge appears with
probability $p$ independently of other edges. By $\vec{B}(n,n,p)$
denote the corresponding space of random directed bipartite graphs
such that each {\bf directed} edge appears with probability $p$
independently of other edges. We also denote by $\{v_1,\ldots,
v_n\}$ and $\{u_1,\ldots, u_n\}$ the corresponding sides.

Now, clearly, if we condition on the event that no edges of the
form $(v_i,u_i)$ or $(u_i,v_i)$, $1 \leq i \leq n$ exist in
$\vec{B}(n,n,p)$ we get a graph distributed precisely as
$H(V,E(H))$. It is also easy to show that the probability of such
an event tends to $e^{-2c}$, thus, it is enough to show that
$\vec{B}(n,n,p)$ has a directed cycle of length
$\Omega(\frac{n}{\ln ^8 n})$ almost surely. We will in fact show
that $B(n,n,p)$ has an {\bf undirected} cycle of that length
almost surely and then use a theorem of McDiarmid \cite{MCD} to
attain the same for $\vec{B}(n,n,p)$. We restate McDiarmid's
result, phrasing it in according to our needs (the theorem is far
more general).

\begin{theo}[McDiarmid]\label{MCDT}
Let $Q$ be a set of directed bipartite graphs with $n$ vertices in
each side that satisfies the following two conditions:
\begin{enumerate}
\item $Q$ is an anti-chain, i.e., for any $D, D' \in Q$, $D \not \subset D'$ and $D' \not
\subset D$.
\item None of the graphs in $Q$ contains opposing edges, i.e. $(u,v)$ and $(v,u)$.
\end{enumerate}
Then,
$$ Pr[\vec{B}(n,n,p) \textrm{ contains a subgraph from $Q$}] \geq
Pr[B(n,n,p) \textrm{ contains a subgraph from $Q$}] .$$
\end{theo}

\noindent It is easy to observe that the set $Q=\{\textrm{all
cycles of length} > \Omega(\frac{n}{\ln ^8 n}) \}$ satisfies the
conditions of Theorem \ref{MCDT}, and therefore, we are left to
prove that $B(n,n,p)$ has an undirected cycle of length
$\Omega(\frac{n}{\ln ^8 n})$.

\begin{prop}\label{PART}
Let $G=(V,E)$ be a connected graph on $m$ vertices with maximum
degree at most $k$. Then for every natural $l$, there exist
disjoint vertex sets $V_1,\ldots,V_t \subset V$, with the
following properties:
\begin{enumerate}
\item $lk \leq |V_i| \leq lk^2$ for every $1 \leq i \leq t$.
\item $\sum ^{t}_{i=1} |V_i| \geq m-lk$.
\item $G[V_i]$ is connected for every $1 \leq i \leq t$.
\end{enumerate}
\end{prop}
\noindent{\bf Proof.} By induction on $m$. For $lk \leq m \leq
lk^2$ we take the whole graph to be $V_1$. For $m > lk^2$, we
choose an arbitrary vertex $v$ and build from it a BFS tree in
which $v$ is the root. For every vertex $w$, denote by $D(w)$ the
set of $w$'s descendants in the BFS tree. We claim that there
exists $w$ such that $lk \leq |D(w)| \leq lk^2$. Otherwise, choose
$w$ such that $|D(w)|\geq lk$ and is minimal, it follows by our
assumption that $|D(w)|>lk^2$. Since the maximum degree is no more
than $k$, $w$ has no more than $k$ direct children. Thus, since
$|D(w)|>lk^2$, one of $w$'s children, $w'$, must have $|D(w')|\geq
lk$ and of course $|D(w')|<|D(w)|$, thus contradicting the
minimality of $|D(w)|$. Now, we take a $w$ such that $lk \leq
|D(w)| \leq lk^2$ and take $V_1 = D(w)$, then $V_1$ is clearly
connected . We remove $D(w)$ from the graph and since the
remaining graph is still connected, we apply the induction
hypothesis to it. $\Box$
\newline

It is a well known fact that if $p=\frac{c}{n}$, where $c>1$, then
$B(n,n,p)$ almost surely has a connected component of size $> Dn$,
(often referred to as the {\em giant component}) where $D=D(c)>0$
is constant. It is also very easy verify that for such $p$, the
maximum degree of $B(n,n,p)$ is almost surely no more than $\ln
n$. Let $H \in B(n,n,p)$ be such a graph, denote its sides by
$A,B$ ($|A|=|B|=n$) and apply Proposition \ref{PART} to its giant
component with $k=\ln n$ and $l=ln ^6 n$. Thus we have $V_1,
\ldots, V_t$, such that $\ln ^7 n \leq V_i \leq \ln ^8 n$,
$H[V_i]$ is connected for every $1 \leq i \leq t$ and $\sum
^{t}_{i=1} |V_i| \geq Dn-\ln ^7 n$. Also, $|V_i \cap A|>\ln ^5 n$
and $|V_i \cap B|>\ln ^5 n$ since otherwise, recalling that the
maximum degree is at most $\ln n$, we would have $|V_i| \leq \ln
^5 n + \ln ^6 n < \ln ^7 n$, for large enough $n$.
\newline

We use now the usual double exposure routine. We first draw
$B(n,n,p=\frac{c'}{n})$ where $1 < c' < c$, then
$B(n,n,p=\frac{\epsilon}{n})$ for $\epsilon = \epsilon (c,c') > 0$
satisfying $(1-\frac{c}{n}) =
(1-\frac{c'}{n})(1-\frac{\epsilon}{n})$. Their union is
distributed precisely $B(n,n,p=\frac{c}{n})$. After the first
draw, we get almost surely such $V_1, \ldots, V_t$ as before.
After the second draw, define a new auxiliary graph $\tilde H$
with $V(\tilde H)=\{V_1,\ldots, V_t\}$ and an edge $\{V_i, V_j\}
\in E(\tilde H)$ for $1 \leq i \neq j \leq t$ iff in the second
draw we have drawn an edge between a vertex in $V_i$ and a vertex
in $V_j$. Clearly, $\tilde H$ is a random graph on no more than
$\frac{Dn}{\ln ^7 n}$ vertices, and since each set $V_i$ has $\ln
^5 n$ vertices on each side, the edge probability, $\tilde p$,
satisfies $1-\tilde p = \left ( 1- \frac{\epsilon}{n} \right
)^{\ln ^{10} n}$ implying $\tilde p > \frac{\epsilon \ln ^{10}
n}{4n} \gg \frac {\ln V(\tilde H)}{V(\tilde H)} \sim \frac{\ln ^9
n}{n}$. But, within this probability range, $\tilde H$ has a
Hamilton cycle almost surely (see, e.g., \cite{BB}). Since $\tilde
H$ has at least $\frac{Dn-\ln ^7 n}{\ln ^8 n}$ vertices, and since
$H[V_i]$ is connected for every $1 \leq i \leq t$ this translates
to a cycle in the original graph, $H$, of length at least
$\frac{Dn - \ln ^7 n}{\ln ^8 n}$, thus concluding our proof.
\hfill $\Box$
\newline


%
%

\section {Case $k\geq 3$}

In \cite{AP}, Achlioptas and Peres use a clever refinement of the
second moment method to estimate the probability that there exists
a satisfying assignment to a random $k$-CNF formula.

We imitate Achlioptas and Peres' proof to prove Theorem \ref{th3}.
In fact, although our proof is only slightly different, we have
not found a way to use Achlioptas and Peres' result to obtain
Theorem \ref{th3}. For the sake of completeness, we repeat the
proof here.
\newline


\noindent{\bf Proof of Theorem \ref{th3}.}

To prove the first part of Theorem \ref{th3}, let $X$ be the r.v.
counting the number of disjoint covers of $H=H_1 \cup H_2$, each
drawn from $H_k(n,m=rn)$ for some constant $r>0$. For a fixed
partition $A\uplus B = [n]$, where $|A|=a$ and $|B|=b=n-a$, the
probability that these are disjoint covers is clearly:

$$\left(1-\frac{\binom{a}{k}}{\binom{n}{k}}\right)^{rn}\left(1-\frac{\binom{b}{k}}{\binom{n}{k}}\right)^{rn}.$$

\noindent As $\binom{y}{k}$ is convex, this expression is
maximized when $a=b=\left \lfloor \frac{n}{2} \right \rfloor$,
which is then approximately, $\left
(1-\frac{1}{2^k}\right)^{2rn}$. Thus,

$$ E[X] \leq \sum _{A\uplus B = [n]} Pr[A\uplus B \textrm{ are disjoint covers}]
\leq \left [2 \left(1-\frac{1}{2^k}\right)^{2r} \right]^n ;$$

\noindent So, if $2 \left(1-\frac{1}{2^k}\right)^{2r} <1$ then
$E[X]=o(1)$ and almost surely, no disjoint covers exist. This
implies that for $r>2^{k-1}\ln 2$, no disjoint covers exist almost
surely.
\newline

Now we prove the second part of the Theorem. We will assume during
the proof that $n$ is a large enough even integer. To get the
theorem for all sufficiently large $n$, observe that if for $r^*$
and even $n$, $H_1, H_2 \sim H_k(n,r^*n)$ have almost surely
disjoint covers, then for any $r<r^*$ and odd $n$, $H_1, H_2 \sim
H_k(n,rn)$ have almost surely disjoint covers. This is because if
we first draw $H_1, H_2 \sim H_k(n+1,r^*(n+1))$ and then delete an
arbitrary vertex we get a hypergraph on $n$ vertices, having
almost surely $r^*n-o(n)$ random edges that almost surely admits
disjoint covers.

A partition $A\uplus B=[n]$ is called {\em balanced} if
$|A|=|B|=n/2$. For a hypergraph $H=H_1 \cup H_2$ denote by $D(H)$
the set of all balanced disjoint covers of $H$. Similarly, for a
blue (red) edge $e$, denote by $D(e)$ the set of all balanced
partitions which colour some vertex of $e$ with blue (red). For a
blue edge, $e$, and a partition, $\sigma$, define $W(\sigma, e)$
to be the number of blue vertices in $e$ minus the number of red
vertices in $e$, with respect to $\sigma$. Similarly, define
$W(\sigma,e)$ for a red edge, $e$, as the number of red vertices
in $e$ minus the number of blue vertices in $e$, with respect to
$\sigma$. For a blue (red) edge, $e$ and a vertex $v\in e$, define
$W(\sigma,e,v)$ to be $1$ if $v$ is coloured blue (red) by
$\sigma$ and $-1$ otherwise. Note that $W(\sigma,e)$ is an integer
between $-k$ and $k$, and if $e=(v_1,\ldots, v_k)$, then
$W(\sigma, e)=\sum_{i=1}^k W(\sigma,e,v_i)$. For constant $0 <
\gamma \leq 1$ to be determined later define the following random
variable:

$$ X = \sum _{\sigma:A\uplus B=[n]} \left(\prod_{e
\textrm{ red}} \gamma^{W(\sigma, e)}\right) \left(\prod_{e
\textrm{ blue}} \gamma^{W(\sigma, e)}\right){\bf 1}_{\{\sigma \in
D(H) \}};$$

\noindent Note that $X>0$ implies $S(D) \neq \emptyset$, i.e.
disjoint covers exist. We will use the following inequality, which
easily follows from Chebyschev's inequality: for any non-negative
random variable $X$,

$$ Pr[X>0] \geq \frac{E[X]^2}{E[X^2]}.$$

Thus, to prove Theorem \ref{th3} it is enough to prove that
$\frac{E[X]^2}{E[X^2]} > C$ for some fixed $C$. To simplify the
analysis, we draw an edge by drawing uniformly at random $k$
vertices with replacements (it is easy to see that the number of
"defect" edges, i.e., edges of size smaller than $k$ is $o(n)$ and
hence the same argument used for dealing with odd $n$ will work
here as well). We begin by evaluating the first moment.
\newline

For any random edge $e$ (either blue or red) and fixed balanced
partition $\sigma$,

$$ E \left [\gamma^{W(\sigma, e)}{\bf 1}_{\{\sigma \in D(e)\}} \right ] =
E \left [\gamma^{W(\sigma,e)} \right ] - \gamma^{-k}Pr[{\bf
1}_{\{\sigma \not \in D(e)\}}] = \left (\frac{\gamma +
\gamma^{-1}}{2}\right)^k - (2\gamma)^{-k}.$$

\noindent Preserving the notation in \cite{AP}, define $\psi
(\gamma) := \left (\frac{\gamma + \gamma^{-1}}{2}\right)^k -
(2\gamma)^{-k}$. Since edges are drawn independently, it follows
that:

$$E[X] = \sum _{\sigma:A\uplus B=[n]} E\left [\left(\prod_{e \textrm{
red}} \gamma^{W(\sigma, e)}\right) \left(\prod_{e \textrm{ blue}}
\gamma^{W(\sigma, e)}\right){\bf 1}_{\{\sigma \in D(H) \}} \right
] = 2\binom{n}{n/2} \psi(\gamma)^{2rn}.$$

\noindent We now evaluate the second moment. Let $\sigma, \tau$ be
two fixed, balanced partitions which have in common precisely
$\alpha n$ vertices in the first part, and let
$e=(v_1,\ldots,v_k)$ be a blue edge. Note that since $\sigma,
\tau$ are balanced, both also have in common precisely $\alpha n$
vertices in the second part and $(1-\alpha )n$ vertices on which
they do not "agree" upon. Then,

$$E\left [\gamma ^{W(\sigma, e)+W(\tau, e)} \right] = \prod
^{k}_{i=1} E\left [\gamma ^{W(\sigma, e, v_i)+W(\tau, e, v_i)}
\right ] = \left [1-2\alpha+\alpha(\gamma^2 + \gamma^{-2}) \right
]^k ,$$

$$ E\left [\gamma ^{W(\sigma, e)+W(\tau, e)}{\bf 1}_{\{\sigma \not
\in D(e)\}} \right ] = \prod ^{k}_{i=1} E \left [ \gamma
^{W(\sigma, e, v_i)+W(\tau, e, v_i)}{\bf 1}_{\{\sigma \textrm{
colours $v_i$ red}\}} \right ] = \left [ \alpha \gamma^{-2} +
\frac{1-2\alpha}{2} \right ]^k ,$$

$$ E \left [ \gamma ^{W(\sigma, e)+W(\tau, e)}{\bf 1}_{\{\sigma ,
\tau \not \in D(e) \}} \right ] = \prod ^{k}_{i=1} E \left [
\gamma ^{W(\sigma, e, v_i)+W(\tau, e, v_i)}{\bf 1}_{\{\sigma ,
\tau \textrm {colour $v_i$ red}\}} \right ] =
\alpha^k\gamma^{-2k}.$$

\noindent So,

\begin{eqnarray*}
E \left [ \gamma ^{W(\sigma, e)+W(\tau, e)}{\bf 1}_{\{\sigma ,
\tau \in D(e) \}} \right ] &=& E \left [ \gamma ^{W(\sigma,
e)+W(\tau, e)}( 1 - {\bf 1}_{\{\sigma \not \in D(e)\}} - {\bf
1}_{\{\tau \not \in D(e)\}} + {\bf 1}_{\{\sigma , \tau \not
\in D(e) \}}) \right ] \\
&=& \left [1-2\alpha+\alpha(\gamma^2 + \gamma^{-2}) \right ]^k -
2\left [ \alpha \gamma^{-2} + \frac{1-2\alpha}{2} \right ]^k +
\alpha^k\gamma^{-2k}.
\end{eqnarray*}

\noindent Denote $f(\alpha) = \left [1-2\alpha+\alpha(\gamma^2 +
\gamma^{-2}) \right ]^k - 2\left [ \alpha \gamma^{-2} +
\frac{1-2\alpha}{2} \right ]^k + \alpha^k\gamma^{-2k}$, then,

\begin{eqnarray*}
E[X^2] &=& \sum _{\sigma , \tau} E \left[ \left ( \prod_{e
\textrm{ red}} \gamma ^{W(\sigma, e)+W(\tau, e)}{\bf 1}_{\{\sigma
, \tau \in D(e) \}}\right ) \left ( \prod_{e \textrm{ blue}}
\gamma ^{W(\sigma, e)+W(\tau, e)}{\bf 1}_{\{\sigma , \tau \in D(e)
\}} \right ) \right ] \\
&=& 2\binom{n}{n/2} \sum ^{n/2}_{z=0} \binom{n/2}{z}^2
f(z/n)^{2rn} .
\end{eqnarray*}

\noindent To estimate this sum we use a lemma proved in \cite{AM}
which is based on the Laplace method for asymptotic integrals
\cite{deB}.

\begin{lemma}\label{LAPLACE}
Let $\phi$ be a positive, twice differentiable function on $[0,1]$
and let $q \geq 1$ be a fixed integer. Let,

$$ S_n = \sum ^{n/q} _{z=0} \binom{n/q}{z} ^q \phi(zq/n)^n .$$

Letting $0^0 := 1$, define $g$ on $[0,1]$ as

$$ g(\alpha)= \frac {\phi(\alpha)}{\alpha^\alpha(1-\alpha)^{1-\alpha}} .$$

If there exists $\alpha _{max} \in (0,1)$ such that $g(\alpha
_{max}) := g_{max} > g(\alpha)$ for all $\alpha \neq \alpha
_{max}$, and $g''(\alpha _{max})<0$, then there exists a constant
$C=C(q, g_{max}, g''(\alpha _{max}),\alpha _{max}) > 0$ such that
$$ S_n < Cn^{-(q-1)/2}g_{max}^n .$$
\end{lemma}

\noindent We apply Lemma \ref{LAPLACE} with $q=2$, $\phi(\alpha) =
f(\alpha /2)^{2r}$. A considerable piece of \cite{AP} is devoted
to prove that for all $k \geq 3$, if $r<2^{k-1}\ln 2 - (\ln
2)(k+1) - \frac{1}{2} - \frac {3}{2k}$, there exists $\gamma$ such
that $\phi(\alpha)$ satisfies the conditions of Lemma
\ref{LAPLACE} with $\alpha _{max} = 1/2$ (see Lemma $3$ of
\cite{AP}). Note that since $f(1/4)=\psi(\gamma)^2$, it follows by
Lemma \ref{LAPLACE} that,

$$ \frac{E[X^2]}{E[X]^2} \leq \frac{2\binom{n}{n/2}Cn^{-1/2}2^n\psi(\gamma)^{4rn}}
{4\binom{n}{n/2}^2\psi(\gamma)^{4rn}} \to \frac
{C}{2}\sqrt{\frac{\pi}{2}},$$

\noindent thus concluding our proof. \hfill $\Box$
\newline

\end{document}